



\magnification=1200
\baselineskip=15pt 


\catcode`@=11 


\font\ninerm=cmr9
\font\eightrm=cmr8
\font\sixrm=cmr6

\font\ninei=cmmi9
\font\eighti=cmmi8
\font\sixi=cmmi6
\skewchar\ninei='177 \skewchar\eighti='177 \skewchar\sixi='177

\font\ninesy=cmsy9
\font\eightsy=cmsy8
\font\sixsy=cmsy6
\skewchar\ninesy='60 \skewchar\eightsy='60 \skewchar\sixsy='60

\font\eightss=cmssq8

\font\ninebf=cmbx9
\font\eightbf=cmbx8
\font\sixbf=cmbx6

\font\ninett=cmtt9
\font\eighttt=cmtt8

\font\ninesl=cmsl9
\font\eightsl=cmsl8

\font\nineit=cmti9
\font\eightit=cmti8


\newskip\ttglue
\def\tenpoint{\def\rm{\fam0\tenrm}%
  \textfont0=\tenrm \scriptfont0=\sevenrm \scriptscriptfont0=\fiverm
  \textfont1=\teni \scriptfont1=\seveni \scriptscriptfont1=\fivei
  \textfont2=\tensy \scriptfont2=\sevensy \scriptscriptfont2=\fivesy
  \textfont3=\tenex \scriptfont3=\tenex \scriptscriptfont3=\tenex
  \def\it{\fam\itfam\tenit}%
  \textfont\itfam=\tenit
  \def\sl{\fam\slfam\tensl}%
  \textfont\slfam=\tensl
  \def\bf{\fam\bffam\tenbf}%
  \textfont\bffam=\tenbf \scriptfont\bffam=\sevenbf
   \scriptscriptfont\bffam=\fivebf
  \def\tt{\fam\ttfam\tentt}%
  \textfont\ttfam=\tentt
  \tt \ttglue=.5em plus.25em minus.15em

  \normalbaselineskip=12pt
  \def\MF{{\manual META}\-{\manual FONT}}%
  \let\sc=\eightrm
  \let\big=\tenbig
  \setbox\strutbox=\hbox{\vrule height8.5pt depth3.5pt width\z@}%
  \normalbaselines\rm}

\def\ninepoint{\def\rm{\fam0\ninerm}%
  \textfont0=\ninerm \scriptfont0=\sixrm \scriptscriptfont0=\fiverm
  \textfont1=\ninei \scriptfont1=\sixi \scriptscriptfont1=\fivei
  \textfont2=\ninesy \scriptfont2=\sixsy \scriptscriptfont2=\fivesy
  \textfont3=\tenex \scriptfont3=\tenex \scriptscriptfont3=\tenex
  \def\it{\fam\itfam\nineit}%
  \textfont\itfam=\nineit
  \def\sl{\fam\slfam\ninesl}%
  \textfont\slfam=\ninesl
  \def\bf{\fam\bffam\ninebf}%
  \textfont\bffam=\ninebf \scriptfont\bffam=\sixbf
   \scriptscriptfont\bffam=\fivebf
  \def\tt{\fam\ttfam\ninett}%
  \textfont\ttfam=\ninett
  \tt \ttglue=.5em plus.25em minus.15em
  \normalbaselineskip=11pt
  \def\MF{{\manual hijk}\-{\manual lmnj}}%
  \let\sc=\ninesc
  \let\big=\ninebig
  \setbox\strutbox=\hbox{\vrule height8pt depth3pt width\z@}%
  \normalbaselines\rm}

\def\eightpoint{\def\rm{\fam0\eightrm}%
  \textfont0=\eightrm \scriptfont0=\sixrm \scriptscriptfont0=\fiverm
  \textfont1=\eighti \scriptfont1=\sixi \scriptscriptfont1=\fivei
  \textfont2=\eightsy \scriptfont2=\sixsy \scriptscriptfont2=\fivesy
  \textfont3=\tenex \scriptfont3=\tenex \scriptscriptfont3=\tenex
  \def\it{\fam\itfam\eightit}%
  \textfont\itfam=\eightit
  \def\sl{\fam\slfam\eightsl}%
  \textfont\slfam=\eightsl
  \def\bf{\fam\bffam\eightbf}%
  \textfont\bffam=\eightbf \scriptfont\bffam=\sixbf
   \scriptscriptfont\bffam=\fivebf
  \def\tt{\fam\ttfam\eighttt}%
  \textfont\ttfam=\eighttt
  \tt \ttglue=.5em plus.25em minus.15em
  \normalbaselineskip=9pt
  \def\MF{{\manual opqr}\-{\manual stuq}}%
  \let\sc=\eightsc
  \let\big=\eightbig
  \setbox\strutbox=\hbox{\vrule height7pt depth2pt width\z@}%
  \normalbaselines\rm}

\def\tenbig#1{{\hbox{$\left#1\vbox to8.5pt{}\right.\n@space$}}}
\def\ninebig#1{{\hbox{$\textfont0=\tenrm\textfont2=\tensy
  \left#1\vbox to7.25pt{}\right.\n@space$}}}
\def\eightbig#1{{\hbox{$\textfont0=\ninerm\textfont2=\ninesy
  \left#1\vbox to6.5pt{}\right.\n@space$}}}


\def\fourteenpointbold{\def\rm{\fam0\frtnb}%
  \textfont0=\frtnb \scriptfont0=\tenbf \scriptscriptfont0=\sevenbf
  \textfont1=\frtnbmi \scriptfont1=\tenbmi \scriptscriptfont1=\sevenbmi
  \textfont2=\frtnbsy \scriptfont2=\tenbsy \scriptscriptfont2=\sevenbsy
  \textfont3=\tenex \scriptfont3=\tenex \scriptscriptfont3=\tenex
  \def\it{\fam\itfam\frtnbit}%
  \textfont\itfam=\frtnbit
  \def\sl{\fam\slfam\frtnbsl}%
  \textfont\slfam=\frtnbsl
  \def\bf{\fam\bffam\frtnrm}
  \textfont\bffam=\frtnrm \scriptfont\bffam=\tenrm
   \scriptscriptfont\bffam=\sevenrm
  \def\tt{\fam\ttfam\frtntt}%
  \textfont\ttfam=\frtntt
  \tt \ttglue=.5em plus.25em minus.15em

  \normalbaselineskip=17pt
  \def\MF{{\manual META}\-{\manual FONT}}%
  \let\sc=\frtnsc
  \let\big=\tenbig
  \setbox\strutbox=\hbox{\vrule height8.5pt depth3.5pt width\z@}%
  \normalbaselines\rm}


\font\frtnb=cmbx12 scaled\magstep1
\font\frtnbmi=cmmib10 scaled\magstep2
\font\frtnbit=cmbxti10 scaled\magstep2
\font\frtnbsl=cmbxsl10 scaled\magstep2
\font\frtnbsy=cmbsy9 scaled\magstep2
\font\frtnrm=cmr12 scaled\magstep1
\font\frtntt=cmtt10 scaled\magstep2
\font\frtnsc=cmcsc10 scaled\magstep2


\font\tenbmi=cmmib10
\font\sevenbmi=cmmib7
\font\tenbsy=cmbsy9
\font\sevenbsy=cmbsy7


\def\ftnote#1{\edef\@sf{\spacefactor\the\spacefactor}#1\@sf
      \insert\footins\bgroup\eightpoint
      \interlinepenalty100 \let\par=\endgraf
        \leftskip=\z@skip \rightskip=\z@skip
        \splittopskip=10pt plus 1pt minus 1pt \floatingpenalty=20000
        \smallskip\item{#1}\bgroup\strut\aftergroup\@foot\let\next}
\skip\footins=12pt plus 2pt minus 4pt 
\dimen\footins=30pc 


\newcount\ftnoteno
\def\Fnote{\global\advance\ftnoteno by 1 \ftnote{$^{\number\ftnoteno}$}}


\input amssym.def \input amssym.tex         


\font\sc=cmcsc10
\font\ninesc=cmcsc9
\font\eightsc=cmcsc8




.263   


\def\head#1#2{\headline{\eightss {\ifnum\pageno=1 \underbar{\raise2pt
 \line {#1 \hfill #2}}\else\hfil \fi}}}

\def\title#1{\centerline{\fourteenpointbold #1}}

\def\author#1{\bigskip\centerline{\sc by #1}}

\def\abstract#1{\vskip.6in\begingroup\ninepoint\narrower\narrower
 \noindent{\bf Abstract.} #1\bigskip\endgroup}

\def\bottom#1#2#3{{\eightpoint\parindent=0pt\parskip=2pt\footnote{}
 {{\it 1991 Mathematics Subject Classification.}
 #1 }\footnote{}{{\it Key words and phrases.} #2 }\footnote{}{#3 }}}

\outer\def\proclaim #1. #2\par{\medbreak \noindent {\sc #1.\enspace }
 \begingroup\it #2 \endgroup
 \par \ifdim \lastskip <\medskipamount \removelastskip \penalty
 55\medskip \fi}  

\def\section#1. #2.{\vskip0pt plus.3\vsize \penalty -150 \vskip0pt
 plus-.3\vsize \bigskip\bigskip \vskip \parskip
 \centerline {\bf \S#1. #2.}\nobreak \medskip \noindent}



\def\proof{\medbreak\noindent{\it Proof.\enspace}}
\def\qed{\quad \vrule height7.5pt width4.17pt depth0pt}
\def\Qed{\qed\ifmmode \relax \else \medbreak \fi}
\long\def\remark{\noindent{\sc Remark.\enspace}}
\long\def\Remark#1\par{\medbreak\remark #1 \medbreak}
\def\example{\noindent{\sc Example:\enspace}}
\def\Example#1\par{\medbreak\example #1 \medbreak}


\def\ref#1 (#2){\par\hangindent=.8cm\hangafter=1\noindent {\sc #1}\ (#2).}
\def\rfr#1 (#2){\par\hangindent=.8cm\hangafter=1\noindent {\sc #1}\ (#2).}
\def\and{{\rm and }}
\def\endreferences{\bigskip\bigskip\endgroup}


\def\beginreferences{\vskip0pt plus.3\vsize \penalty -150 \vskip0pt
 plus-.3\vsize \bigskip\bigskip \vskip \parskip
 \begingroup\baselineskip=12pt\frenchspacing
   \centerline{\bf REFERENCES}
   \vskip12pt\parindent=0pt
   \def\and{{\rm and }}
   \def\rrf##1 (##2){{\sc ##1}\ (##2).}
   \everypar={\hangindent=.8cm\hangafter=1\rrf }}
\def\endreferences{\bigskip\bigskip\endgroup}


\def\implies{\quad\Longrightarrow\quad}
\def\cbuldot{{\,\raise.25ex\hbox{$\scriptscriptstyle\bullet$}}}
\def\I#1{{\bf 1}_{#1}}
\def\E{{\bf E}}

\def\P{{\bf P}}

\def\R{{\bf R}}


\def\la#1{\mathop{#1}\limits^{\lower 6pt\hbox{$\leftarrow$}}}
\def\ra#1{\mathop{#1}\limits^{\lower 6pt\hbox{$\rightarrow$}}}


\overfullrule=0pt

\baselineskip=14.5pt  



\newif\iflefteqnumbers
\lefteqnumbersfalse

\newif\ifproofmode          
\proofmodefalse    

\newif\ifforwardreference      
\forwardreferencefalse            

\newif\ifcheckdefinedreference     
\checkdefinedreferencefalse        

\newif\ifreportundefref           
\reportundefreffalse              

\newif\iflblnonempty         
\lblnonemptyfalse            

%
%

\newif\ifmultilinedisplay  
\multilinedisplayfalse

\newif\ifequationparen     
\equationparenfalse        

\newif\ifchaprefwordy       
\chaprefwordyfalse          

\newif\ifchapnumroman       
\chapnumromanfalse          

\newif\ifsectionnumbers      
\sectionnumbersfalse              

\newif\ifcontinuousnumbering    
\continuousnumberingfalse         

\newif\iffiguresectionnumbers  
\figuresectionnumbersfalse          

\newif\ifcontinuousfigurenumbering   
\continuousfigurenumberingfalse

\newif\iftheoremcounting
\theoremcountingfalse


\font\marginstyle=cmss8      

\font\headfont=cmss8 
\font\sc=cmcsc10

\def\proctitlefont{\sc}    
\def\procbodyfont{\it}     
\def\procnumberfont{\rm}   

\newtoks\sectionletter    
\newcount\sectionno          
\newcount\eqlabelno    
\newcount\figureno      
\newcount\theoremlabelno 
\newcount\chapno         

%
%
%

\def\getlabel#1{\csname#1\endcsname}
\def\ifundefined#1{\expandafter\ifx\csname#1\endcsname\relax}
\def\stripchap#1/#2.#3?{#1}
\def\stripsec#1/#2.#3?{#2}
\def\stripeq#1/#2.#3?{#3}

\sectionno=0
\eqlabelno=0
\figureno=0
\theoremlabelno=0
\chapno=0

%
%

\def\chapfolio{\ifchapnumroman\uppercase\expandafter{\romannumeral\the\chapno}%
                \else \the\chapno\fi}

\def\bumpchapno{\global\advance\chapno by 1 \global\sectionno=0
                \global\eqlabelno=0 \global\figureno=0
                \global\theoremlabelno=0}

%
%

\def\chapshow#1{{\chapno=#1\chapfolio}}  

%
%

\def\sectionfolio{\ifnum \sectionno>0 \the \sectionno \else
                           \the\sectionletter \fi}

%
%

\def\bumpsectionno{\ifnum\sectionno>-1 \global\advance\sectionno by 1
        \else\global\advance\sectionno by -1 \setletter\sectionno\fi
        \ifcontinuousnumbering\else\global\eqlabelno =0
                                 \global\theoremlabelno=0 \fi
        \ifcontinuousfigurenumbering\else\global\figureno=0 \fi}

%
%

\def\setletter#1{\ifcase-#1 {}\or{}    \or\global\sectionletter={A}
      \or\global\sectionletter={B} \or\global\sectionletter={C}
      \or\global\sectionletter={D} \or\global\sectionletter={E}
      \or\global\sectionletter={F} \or\global\sectionletter={G}
      \or\global\sectionletter={H} \or\global\sectionletter={I}
      \or\global\sectionletter={J} \or\global\sectionletter={K}
      \or\global\sectionletter={L} \or\global\sectionletter={M}
      \or\global\sectionletter={N} \or\global\sectionletter={O}
      \or\global\sectionletter={P} \or\global\sectionletter={Q}
      \or\global\sectionletter={R} \or\global\sectionletter={S}
      \or\global\sectionletter={T} \or\global\sectionletter={U}
      \or\global\sectionletter={V} \or\global\sectionletter={W}
      \or\global\sectionletter={X} \or\global\sectionletter={Y}
      \or\global\sectionletter={Z} \fi}

%
%

\def\tempsetletter#1{\ifcase-#1 {}\or{} \or\sectionletter={A}
      \or\sectionletter={B} \or\sectionletter={C}
      \or\sectionletter={D} \or\sectionletter={E}
      \or\sectionletter={F} \or\sectionletter={G}
      \or\sectionletter={H} \or\sectionletter={I}
      \or\sectionletter={J} \or\sectionletter={K}
      \or\sectionletter={L} \or\sectionletter={M}
      \or\sectionletter={N} \or\sectionletter={O}
      \or\sectionletter={P} \or\sectionletter={Q}
      \or\sectionletter={R} \or\sectionletter={S}
      \or\sectionletter={T} \or\sectionletter={U}
      \or\sectionletter={V} \or\sectionletter={W}
      \or\sectionletter={X} \or\sectionletter={Y}
      \or\sectionletter={Z} \fi}

%
%

\def\sectionshow#1{\ifnum#1>0 #1
        \else{\tempsetletter{\number#1}\the\sectionletter}\fi}

%
%

\def\chapsecshow#1/#2{\ifchaprefwordy\sectionshow{#2} of
                    Chapter~\chapshow{#1}%
          \else\chapshow{#1}--\sectionshow{#2}\fi}
%
%
%
%

\def\chapsecitemshow#1/#2.#3{\ifchaprefwordy\ifequationparen
                    (\fi\sectionshow{#2}.#3\ifequationparen)\fi\ of
                    Chapter~\chapshow{#1}%
          \else\ifequationparen(\fi\chapshow{#1}--\sectionshow{#2}.#3%
                 \ifequationparen)\fi\fi}

\def\secitemshow#1.#2{\ifequationparen(\fi\sectionshow{#1}.#2%
                         \ifequationparen)\fi}

\def\chapitemshow#1/#2{\ifchaprefwordy \ifequationparen(\fi#2%
              \ifequationparen)\fi\ of Chapter~\chapshow{#1}%
          \else\ifequationparen(\fi\chapshow{#1}--#2\ifequationparen)\fi\fi}

\def\itemshow#1{\ifequationparen(\fi#1\ifequationparen)\fi}

%
%
\def\today{\number\day\ \ifcase\month\or
    January\or February\or March\or April\or May\or June\or
    July\or August\or September\or October\or November\or December\fi
    \ \number\year}

%
%
\def\initialeqmacro{\ifproofmode  
              \headline{\headfont\underbar{\raise2pt
              \line{\today\hskip 10pt [\the\time]\hfill\jobname\ --- draft}}}
                      \fi
          \openin1=\jobname.lbl
          \ifeof1 \lblnonemptyfalse   
           \else\lblnonemptytrue
            \fi
          \ifforwardreference 
             \iflblnonempty
                   \input \jobname.lbl
                   \fi
             \immediate\openout1=\jobname.lbl
                  \fi
          \ifcheckdefinedreference\iflblnonempty\reportundefreftrue
                                   \else\reportundefreffalse
                                    \fi
           \else\ifforwardreference\reportundefreffalse
                 \fi
             \fi}


\let\END=\end
\def\end{\ifforwardreference\iflblnonempty
                   \immediate\write16{ *** If you have made
                              changes to the labels, you may need}
                   \immediate\write16{     to rerun TeX to get the
                              crossreferences correct. *** }
                 \else\immediate\write16{ *** Rerun TeX
                              to get the crossreferences correct. *** }
                   \fi
           \fi \END}


%
%

\def\chaplabel#1{\bumpchapno\ifforwardreference
      \immediate\write1{\noexpand\expandafter\noexpand\def
      \noexpand\csname CHAPLABEL#1\endcsname{\the\chapno}}\fi
      \global\expandafter\edef\csname CHAPLABEL#1\endcsname
      {\the\chapno}}

\def\chapref#1{\ifundefined{CHAPLABEL#1}Chapter $$
             \ifforwardreference \relax
                  \else \write16{ ***Undefined Chapter Reference #1*** } \fi
          \else\edef\LABxx{\getlabel{CHAPLABEL#1}}
            Chapter~\chapshow\LABxx\fi}

%
%

\def\sectionlabel#1{\bumpsectionno\ifforwardreference
      \immediate\write1{\noexpand\expandafter\noexpand\def
      \noexpand\csname SECLABEL#1\endcsname{\the\chapno/\the\sectionno}}\fi
      \global\expandafter\edef\csname SECLABEL#1\endcsname
      {\the\chapno/\the\sectionno}}

\def\sectionref#1{\ifundefined{SECLABEL#1}Section ??
             \ifforwardreference \relax
                  \else \write16{ ***Undefined Section Reference #1*** } \fi
          \else\edef\LABxx{\getlabel{SECLABEL#1}}
           \edef\chapnoref{\expandafter\stripchap\LABxx.?}%
           \edef\sectionnoref{\expandafter\stripsec\LABxx.?}%
            \ifnum\sectionnoref>-1 Section~\else Appendix\fi%
            \ifnum\chapnoref=\chapno \sectionshow\sectionnoref%
              \else\chapsecshow\chapnoref/\sectionnoref\fi\fi}

%
%

\def\itemref#1{\ifundefined{EQLABEL#1}--?--
         \ifreportundefref
          \write16{ ***Undefined Equation Reference #1*** }
           \fi
      \else \edef\LABcse{\getlabel{EQLABEL#1}}%
       \edef\LABch{\expandafter\stripchap\LABcse}%
       \edef\LABse{\expandafter\stripsec\LABcse}%
       \edef\LABeq{\expandafter\stripeq\LABcse}%
         \ifnum\LABch=\chapno
            \ifsectionnumbers\secitemshow\LABse.\LABeq%
             \else\ifnum \LABse=\sectionno \itemshow\LABeq%
                   \else\ifcontinuousnumbering\itemshow\LABeq%
                         \else\secitemshow\LABse.\LABeq%
                          \fi
                     \fi
                \fi
           \else\ifsectionnumbers\chapsecitemshow\LABch/\LABse.\LABeq%
                 \else\ifcontinuousnumbering\chapitemshow\LABch/\LABeq%
                       \else\chapsecitemshow\LABch/\LABse.\LABeq%
                        \fi
                    \fi
              \fi
        \fi}

%
%

\def\figlabel#1{\global\advance\figureno by 1 \relax
      \ifforwardreference
         \immediate \write1{\noexpand\expandafter\noexpand\def
         \noexpand\csname
         FIGLABEL#1\endcsname{\the\chapno/\the\sectionno.\the\figureno?}}\fi
      \global\expandafter\edef\csname FIGLABEL#1\endcsname
        {\the\chapno/\the\sectionno.\the\figureno?}
      \ifproofmode\llap{\hbox{\marginstyle #1\ }}\relax\fi
      {\bf Figure \iffiguresectionnumbers\sectionfolio.\fi\the\figureno.}}

\def\figref#1{\ifundefined{FIGLABEL#1}!!!!
          \ifreportundefref
           \write16{ ***Undefined Figure Reference #1*** }
            \fi
      \else \edef \LABcsf{\getlabel{FIGLABEL#1}}%
       \edef\LABch{\expandafter\stripchap\LABcsf}%
       \edef\LABse{\expandafter\stripsec\LABcsf}%
       \edef\LABfi{\expandafter\stripeq\LABcsf}%
         \ifnum\LABch=\chapno
            \iffiguresectionnumbers\secitemshow\LABse.\LABfi%
             \else\ifnum \LABse=\sectionno \itemshow\LABfi%
                   \else\ifcontinuousfigurenumbering\itemshow\LABfi%
                         \else\secitemshow\LABse.\LABfi%
                          \fi
                     \fi
                \fi
           \else\ifsectionnumbers\chapsecitemshow\LABch/\LABse.\LABfi%
                 \else\ifcontinuousfigurenumbering\chapitemshow\LABch/\LABfi%
                       \else\chapsecitemshow\LABch/\LABse.\LABfi%
                        \fi
                    \fi
              \fi
        \fi}

%
%

\def\lefteqn#1{\hbox to 0pt{$\displaystyle #1$\hss}}

%
%

%
%

\def\bsection#1#2 {\sectionlabel{#2} \vskip0pt plus.3\vsize \penalty -250
                \vskip0pt plus-.3\vsize \bigskip\bigskip \vskip \parskip
    \centerline{{\bf {}\llap{%
      \ifproofmode{\marginstyle #2\hskip10pt}\else {}\fi}
      \ifnum\sectionno>-1
        \S\the\sectionno.\enspace
       \else Appendix \the\sectionletter:\enspace
      \fi #1.}}\nobreak \smallskip \noindent}%

%
%

\def\createlabel#1{\global\advance\eqlabelno by 1 
      \ifforwardreference\immediate\write1{\noexpand\expandafter
        \noexpand\def \noexpand\csname EQLABEL#1\endcsname{\the
        \chapno/\the\sectionno.\the\eqlabelno?}}\fi
      \global\expandafter\edef\csname EQLABEL#1\endcsname{\the
      \chapno/\the\sectionno.\the\eqlabelno?}}
%

\def\createthmlabel#1{\global\advance\theoremlabelno by 1
      \ifforwardreference\immediate\write1{\noexpand\expandafter
        \noexpand\def \noexpand\csname EQLABEL#1\endcsname{\the
        \chapno/\the\sectionno.\the\theoremlabelno?}}\fi
      \global\expandafter\edef\csname EQLABEL#1\endcsname{\the
      \chapno/\the\sectionno.\the\theoremlabelno?}}


%
%
\def\begineqalno#1\endeqalno{\multilinedisplaytrue
 \iflefteqnumbers $$\leqalignno{#1}$$ 
      \else $$\eqalignno{#1}$$ \fi\multilinedisplayfalse}
\def\endeqalno{}

\def\label #1 {
      \createlabel{#1}
      \ifmultilinedisplay & \fi
         \iflefteqnumbers \ifmultilinedisplay\else\leqno\fi
             \ifproofmode\hskip-.2truein\llap{\marginstyle
                        #1}\hskip .2truein\fi 
             \else \ifmultilinedisplay\else\eqno\fi \fi
      (\ifsectionnumbers\sectionfolio.\fi\the\eqlabelno)
      \ifproofmode \iflefteqnumbers \relax
                       \else \rlap{\marginstyle \hskip .2truein #1 }\fi\fi}


\def\procl#1.#2 {\iftheoremcounting
      \createthmlabel{#1.#2} \else \createlabel{#1.#2}\fi
    \medbreak \noindent
    \ifproofmode
          \hskip-.2truein\llap{\marginstyle #1.#2}\hskip .2truein
    \fi
    \proctitlefont {\ifx#1t\hbox{Theorem}\fi \ifx#1l\hbox{Lemma}\fi
                \ifx#1c\hbox{Corollary}\fi \ifx#1p\hbox{Proposition}\fi
                \ifx#1r\hbox{Remark}\fi \ifx#1d\hbox{Definition}\fi
                \ifx#1x\hbox{Example}\fi \ifx#1g\hbox{Conjecture}\fi}
    \procnumberfont
      \ifsectionnumbers\sectionfolio.\fi
      \iftheoremcounting\the\theoremlabelno.\enspace
                        \else\the\eqlabelno.\enspace\fi
        \ifx #1t  \procbodyfont  \fi
        \ifx #1l  \procbodyfont  \fi
        \ifx #1c  \procbodyfont  \fi
        \ifx #1p  \procbodyfont  \fi 
        \ifx #1g  \procbodyfont  \fi }

\def\endprocl{\rm \par \ifdim \lastskip <\medskipamount \removelastskip
  \penalty 55 \medskip \fi}

%
%

\def\procname#1{\kern-.25em{\proctitlefont (#1)}\enspace}

\def\proof{\medbreak\noindent{\it Proof.\enspace}}

\def\proofof #1.#2 {\medbreak\noindent
     {\it Proof of 
      {\ifx #1tTheorem\fi
        \ifx #1lLemma\fi 
        \ifx #1cCorollary\fi 
        \ifx #1pProposition\fi}
       \itemref {#1.#2}.}\enspace}

%
%

\def\ref#1.#2/{\ifx #1s\sectionref{#1.#2}%
                 \else
                   \ifx #1C\chapref{#1.#2}%
                    \else
                      \ifx #1tTheorem~\fi  
                      \ifx #1lLemma~\fi 
                      \ifx #1cCorollary~\fi 
                      \ifx #1pProposition~\fi 
                      \ifx #1rRemark~\fi 
                      \ifx #1dDefinition~\fi
                      \ifx #1xExample~\fi 
                      \ifx #1fFigure~\fi 
                      \ifx #1e\global\equationparentrue\fi
                     \itemref{#1.#2}\global\equationparenfalse
                    \fi
                 \fi}


\theoremcountingtrue
\continuousfigurenumberingtrue
\forwardreferencetrue
\checkdefinedreferencetrue
\lefteqnumberstrue
\initialeqmacro

 at 8 truept   

\font\addrfont=cmcsc10 at 10 truept

\def\E{{\bf E}}

\def\R{{\bf R}}
\def\P{{\bf P}}
\def\F{{\cal F}}    
\def\G{{\cal G}}    
\def\L{{\cal L}}    
\def\Larrow#1{{\buildrel{\lower1.5pt\hbox{$\scriptscriptstyle\leftarrow$}}
 \over {#1}}}
\def\Rarrow#1{{\buildrel{\lower1.5pt\hbox{$\scriptscriptstyle\rightarrow$}}
 \over {#1}}}
\def\longRarrow#1{{\buildrel{\lower1.5pt\hbox{$\scriptscriptstyle
  \longrightarrow$}} \over {#1}}}

\def\Lh{\widehat \L}
\def\muh{{\widehat{\mu}}}  
\def\muhs{\muh{}^*}  
\def\Fs{\F{}^*}  
\def\l#1{\langle #1, \L \rangle}

\ifproofmode \relax \else\head{Classical and Modern Branching Processes}
{Version of 21 June 1995}\fi

\vglue20pt
\title{A Simple Path to Biggins' Martingale Convergence}
\title{for Branching Random Walk}
\author{Russell Lyons}
\abstract{We give a simple non-analytic proof of Biggins' theorem on
martingale convergence for branching random walks.}

\bottom{Primary 60J80.}{Galton-Watson.}
{Research partially supported by the Institute for Mathematics and Its
Applications (Minneapolis) and NSF Grant DMS-93069.}

Let $\L := \{X_i\}_{i=1}^L$ be a random $L$-tuple of real numbers, where
$L$ is also random and can take the values 0 and $\infty$. This
can also be thought of as an ordered
point process on $\R$.   The random variable $\L$
is used as the basis
for construction of a branching random walk in the usual way: 
An initial particle at the origin of $\R$ gives birth to $L$
particles with displacements $X_1, X_2, \ldots$. Then each of these
particles gives birth to a random number of particles with random
displacements from its new position according to the same law as $\L$ and
independently of one another and of the initial displacements.
This continues in a like manner forever or until there are no more
particles. For a particle $\sigma$, write $|\sigma|$ for the generation
in which $\sigma$ is born, $X(\sigma)$ for its displacement from its
parent, and $S(\sigma)$ for its
position. Denote the initial particle by 0, also known as the root of the
family tree.  If $\tau$ is an ancestor of $\sigma$, write $\tau <
\sigma$. Thus, we have
$$
S(\sigma) = \sum_{0 < \tau \le \sigma} X(\tau) \,.
$$
Also, write $\L(\sigma)$ for the copy of $\L$ used to generate
the children of $\sigma$.
Let $q$ be the extinction probability of the underlying Galton-Watson
process. 

 For $\alpha \in \R$, define $\l{\alpha} := \sum_{i=1}^L
e^{-\alpha X_i}$ and $m(\alpha) := \E[\l{\alpha}] \in (0, 
\infty]$. Assume that $m(0) > 1$, so that $q < 1$.
If $m(\alpha) < \infty$ for some $\alpha$, then the sequence 
$$
W_n(\alpha) := {\sum_{|\sigma| = n} e^{-\alpha S(\sigma)} \over
                m(\alpha)^n}
$$
is a martingale with a.s.\ limit $W(\alpha)$. Write
$$
m'(\alpha) := \E\left[\sum_{i=1}^L X_i e^{-\alpha X_i}\right]
$$
when this exists in $[-\infty,\,\infty]$ as a Lebesgue integral.
Biggins (1977) has determined when $W(\alpha)$ is nontrivial:

\proclaim Biggins' Theorem. Suppose that $\alpha \in \R$ is such that
$m(\alpha) < \infty$ and $m'(\alpha)$ exists and is finite.
Then the following are equivalent:
\smallskip
\item{(i)} $\P[W(\alpha) = 0] = q$;
\item{(ii)} $\P[W(\alpha) = 0] < 1$;
\item{(iii)} $\E[W(\alpha)] = 1$;
\item{(iv)} $\E[\l{\alpha} \log^+ \l{\alpha}] < \infty$
          and ${\alpha m'(\alpha)/m(\alpha)} < \log m(\alpha)$.


\Remark In fact, the hypotheses here are very slightly weaker than those of
Biggins (1977), Lemma 5. Moreover, the proof to follow works without the
assumption that $m'(\alpha)$ be finite, except for the implication
(ii) $\Rightarrow$ (iv), where it needs the assumption that $\alpha m'(\alpha)
\ne -\infty$.

\Remark The case of Biggins' Theorem where $L$ is constant, $X_i$ are 
independent and identically distributed,
and $m(\alpha) = 1$ was proved also by Kahane (see
Kahane and Peyri\`ere (1976); the first condition in (iv) above follows
from the assumptions that $m(\alpha) < \infty$ and $|m'(\alpha)| < \infty$
in the case that $L$ is bounded since convexity of the function $x
\mapsto x \log x$ shows that $\E[\l{\alpha} \log^+ \l{\alpha}] \le
|m'(\alpha)| + \| \log^+ L \|_\infty m(\alpha)$).
When the conclusions of Biggins' Theorem hold in Kahane's context,
the measure $\muh$ below is introduced by Peyri\`ere on p.~141
of Kahane and Peyri\`ere (1976) for another purpose.
It and related constructions in other situations also occur, usually
including the same direct construction as ours, in 
Kallenberg (1977), Hawkes (1981), Rouault (1981),
Joffe and Waugh (1982), Kesten (1986),
Chauvin and Rouault (1988), Chauvin, Rouault and
Wakolbinger (1991), and Waymire and Williams (1995).


Evidently, (iii) implies (ii).
The fact that (i) and (ii) are equivalent follows from the standard
``zero-one" property of Galton-Watson processes. We shall present a simple
proof of the other equivalences modelled on the proof of the
Kesten-Stigum theorem in Lyons, Pemantle and Peres (1995). I am grateful
to Anatole Joffe for asking me for the details of how this is done.
The same method
was discovered independently by Waymire and Williams (1995)
for the case treated by Kahane (mentioned above). In fact, Waymire and
Williams relax the condition that the $X_i$ be i.i.d. They even relax the
independence of the $\L(\sigma)$, which could be done here as well.

\proof Fix $\alpha$.
 If $t$ is a rooted tree (with distinguishable vertices)
 and $X$ is a real-valued function on the vertices of
$t$ other than its root, we call $(t, X)$ a {\bf labelled tree}. A {\bf
ray} in a tree is an infinite
line of descent starting from the root. Given a ray $\xi$,
the vertex on $\xi$ in generation $n$ is denoted $\xi_n$.
 In the space of
labelled trees, let $\F_n$ denote the $\sigma$-field generated by the
first $n$ levels. We shall also work on the space of labelled trees with
distinguished rays, $(t, X, \xi)$;
denote by $\Fs_n$ the $\sigma$-field generated by the
first $n$ levels there.
For $\sigma \in t$, write $S(\sigma) := \sum_{0 < \tau \le \sigma} X(\tau)$
and set
$$
W_n(t, X) := {\sum_{|\sigma| = n} e^{-\alpha S(\sigma)} \over
                     m(\alpha)^n}\,.
$$

Branching random walk gives a probability
measure, $\mu$, on the set of labelled trees.
We shall construct a related probability
measure $\muhs$ on the set of infinite labelled trees with
distinguished rays. Let $\mu_n$ be the restriction of $\mu$ to $\F_n$.
Now any $\Fs_n$-measurable function $f$ can be written as
$$
f(t, X, \xi) = \sum_{|\sigma|=n} f_\sigma(t, X) \I{\xi_n = \sigma}
$$
for some $\F_n$-measurable functions $f_\sigma$. Let
$\mu^*_n$ be counting measure on $\{|\sigma|=n\}$ fibered over $\mu_n$; more
precisely, $\mu^*_n$ is the (non-probability) measure on $\Fs_n$ such that
for all nonnegative $\Fs_n$-measurable functions $f$,
$$
\int f(t, X, \xi) \,d\mu^*_n = \int \sum_{|\sigma|=n} f_\sigma(t, X)
\,d\mu_n \,.
$$
Then the measure $\muhs$ will satisfy
$$
 {d\muhs_n \over d\mu^*_n} (t, X,\xi) =
      {e^{-\alpha S(\xi_n)} \over m(\alpha)^n}
    \label e.goal
$$
 for all $n$ and all labelled trees with rays $(t, X, \xi)$, where
$\muhs_n$ denotes the restriction of $\muhs$ to $\Fs_n$.
The projection of $\muhs$ to the space of trees, denoted
$\muh$, then satisfies
$$
{d\muh_n \over d\mu_n} (t, X) = W_n(t, X) \label e.rn
$$
for all $n$ and all labelled trees $(t, X)$, where $\muh_n$ denotes the
restriction of $\muh$ to $\F_n$.

To define $\muhs$,
let $\Lh$ be a random variable whose law has Radon-Nikodym derivative
$\l{\alpha}/m(\alpha)$ with respect to the law of $\L$.
 Start with one particle $v_0$ at the origin.  Generate
offspring and displacements
according to a copy $\Lh_1$ of $\Lh$.
Pick one of these children $v_1$ at random, where a child is picked with 
probability proportional to $e^{-\alpha X}$ when its displacement is
$X$. The children other than $v_1$
give rise to ordinary independent branching random walks, while $v_1$
gets an independent number of offspring and displacements
according to a copy $\Lh_2$ of $\Lh$.  Again,
pick one of the children
of $v_1$ at random, call it $v_2$, with the others giving rise to
ordinary independent branching random walks, and so on.

Define the measure $\muhs$ as the joint
distribution of the random labelled tree and the random ray
$ (v_0,v_1,v_2,\ldots)$. Then
$$
{d\muh^*_{n+1} \over d\mu^*_{n+1}} (t, X, \xi) 
  = {d\muh^*_{n} \over d\mu^*_{n}} (t, X, \xi) \cdot {\langle \alpha,
\L(\xi_n) \rangle \over m(\alpha)} \cdot {e^{-\alpha X(\xi_{n+1})} \over
\langle \alpha, \L(\xi_n) \rangle}
  = {e^{-\alpha X(\xi_{n+1})} \over m(\alpha)} {d\muh^*_{n} \over
d\mu^*_{n}} (t, X, \xi) \,.
$$
Thus, \ref e.goal/ follows by induction.

Note that for any $k \ge 0$,
$$
\int X(v_k) \,d\muhs = \E\left[\sum_{i=1}^L X_i {e^{-\alpha X_i} \over
\l{\alpha}} {\l{\alpha} \over m(\alpha)} \right] = -m'(\alpha)/m(\alpha)
\,. \label e.mean
$$
Thus, by the strong law of large numbers, we have
$$
S(v_n)/n \to -m'(\alpha)/m(\alpha) \quad \muhs\hbox{-a.s.} \label e.sl
$$

For any labelled tree $(t, X)$, set $W(t, X) := \limsup W_n(t,
X)$. Now by \ref e.rn/, we have the implications (Durrett (1991), p. 210,
Exercise 3.6)
$$
W(t, X) = \infty \quad \muh\hbox{-a.s.} \implies W(t, X) = 0 \quad
\mu\hbox{-a.s.} \,, \label e.1
$$
$$
W(t, X) < \infty \quad \muh\hbox{-a.s.} \implies \int W(t, X) \,d\mu = 1
\,.
\label e.2
$$

Suppose first that (iv) fails. We have
$$
W_{n+1}(t, X) \ge {e^{-\alpha S(v_{n})} \over m(\alpha)^{n+1}} \langle
\alpha, \Lh_{n+1} \rangle \,, \label e.lb
$$
with the two terms in the product being $\muhs$-independent.
Now, if ${\alpha m'(\alpha)/m(\alpha)} \ge \log m(\alpha)$, then 
$\limsup e^{-\alpha S(v_{n})} / m(\alpha)^{n} = \infty$ by \ref e.sl/ in
case ${\alpha m'(\alpha)/m(\alpha)} > \log m(\alpha)$
and by \ref e.mean/ and the Chung-Fuchs theorem in case
${\alpha m'(\alpha)/m(\alpha)} = \log m(\alpha)$.
This implies by \ref e.lb/ that $W(t, X) = \infty$
$\muh$-a.s., whence by \ref e.1/, (ii) fails. On the other hand,
if ${\alpha m'(\alpha)/m(\alpha)} < \log m(\alpha)$, then since $\E[\log^+
\langle \alpha, \Lh \rangle] = \E[\l{\alpha} \log^+ \l{\alpha}]/m(\alpha)
= \infty$ by assumption, we have that 
$$
\limsup {1 \over n} \log^+\langle \alpha, \Lh_n
\rangle  = \infty \quad \muhs\hbox{-a.s.}
$$
by virtue of the Borel-Cantelli lemma.
This means that the first term in the right-hand side of \ref e.lb/ decays 
exponentially while the second has superexponential explosions. Hence,
again, $W(t, X) = \infty$ $\muh$-a.s.

Conversely, suppose that (iv) holds. Let $\G$ be the $\sigma$-field generated
by $\{\Lh_k\}_{k \ge 1}$. Then
\begineqalno
\E_{\muhs}[W_n(t, X) \mid \G] &= \sum_{k=0}^{n-1} {e^{-\alpha S(v_k)} \over
               m(\alpha)^{k+1}} \biggl(\langle \alpha, \Lh_{k+1} \rangle
                - e^{-\alpha X(v_{k+1})} \biggr) +
                     {e^{-\alpha S(v_n)} \over m(\alpha)^n} \cr
  &= \sum_{k=0}^{n-1} {e^{-\alpha S(v_k)} \over m(\alpha)^{k+1}}
                      \langle \alpha, \Lh_{k+1} \rangle
    - \sum_{k=1}^{n-1} {e^{-\alpha S(v_k)} \over m(\alpha)^k} \,.\cr
\endeqalno
By hypothesis and \ref e.sl/,
the terms ${e^{-\alpha S(v_k)} / m(\alpha)^k}$ decay
exponentially while the terms $\langle \alpha, \Lh_{k+1} \rangle$ grow
(at most) subexponentially by the Borel-Cantelli lemma again.
Therefore both series converge $\muhs$-a.s., whence
$\liminf W_n(t, X) < \infty$ $\muh$-a.s. by Fatou's lemma.
In light of \ref e.rn/, $\{1/W_n(t, X)\}$ is a $\muh$-martingale, so that
$\{W_n(t, X)\}$ converges $\muh$-a.s.
Thus, we have $W(t, X) < \infty$ $\muh$-a.s. and (iii) is a
consequence of \ref e.2/. \Qed

\beginreferences

Biggins, J. D. (1977) Martingale convergence in the branching random
walk. {\it J. Appl. Prob.} {\bf 14}, 25--37.

Chauvin, B. \and Rouault, A. (1988) KPP equation and supercritical
branching Brownian motion in the subcritical speed area. Application to
spatial trees, {\it Probab. Theory Relat. Fields} {\bf 80}, 299--314.

Chauvin, B., Rouault, A., \and Wakolbinger, A. (1991) Growing
conditioned trees, {\it Stochastic Process. Appl.} {\bf 39}, 117--130.

Durrett, R. (1991) Probability: Theory and Examples. Wadsworth,
Pacific Grove, California.

Hawkes, J. (1981) Trees generated by a simple branching process, {\it
J. London Math. Soc.} {\bf 24}, 373--384.

Joffe, A. \and Waugh, W. A. O'N. (1982) Exact distributions of kin
numbers in a Galton-Watson process, {\it J. Appl. Prob.} {\bf 19},
767--775.

Kahane, J.-P. \and Peyri\`ere, J. (1976) Sur certaines martingales de
Benoit Mandelbrot, {\it Adv. in Math.} {\bf 22}, 131--145.

Kallenberg, O. (1977) Stability of critical cluster fields, {\it Math.
Nachr.} {\bf 77}, 7--43.

Kesten, H. (1986) Subdiffusive behavior of random walk on
a random cluster, {\it Ann. Inst. H. Poincar\'e Probab. Statist.}
{\bf 22}, 425--487.

Lyons, R., Pemantle, R. \and Peres, Y. (1995) Conceptual proofs of $L
\log L$ criteria for mean behavior of branching processes, {\it Ann.
Probab.}, to appear.

Rouault, A. (1981) Lois empiriques dans les processus de branchement
spatiaux homog\`enes supercritiques, {\it C. R. Acad. Sci. Paris S\'er.
I. Math.} {\bf 292}, 933--936.

Waymire, E. C. \and Williams, S. C. (1995) A cascade decomposition theory
with applications to Markov and exchangeable cascades, {\it Trans. Amer.
Math. Soc.}, to appear.

\endreferences
\nobreak
\vskip-10pt
\begingroup
\addrfont
\parindent=0pt\baselineskip=10pt

\hfill
Department of Mathematics,
Indiana University,
Bloomington, IN 47405-5701

\endgroup

\bye